\documentclass[12pt]{article} \textwidth=6.2in \oddsidemargin=0in
\begin{document} 
\def\ov{\over} \def\t{\tau} \def\s{\sigma} \def\sp{\vspace{1ex}}
\def\cd{\cdots} \def\iy{\infty}\def\inv{^{-1}} \def\dl{\delta}
\def\be{\begin{equation}} \def\ee{\end{equation}} \def\s{\sigma} 
\def\ph{\varphi} \def\ps{\psi} \def\pl{\partial}\def\r{\rho} 
\def\ep{\varepsilon} \def\sech{{\rm sech}\,} \def\I{{\rm Im}\,}
\def\R{{\rm Re}\,} \def\ep{\varepsilon} \def\a{\alpha} \def\G{\Gamma}
\def \tlG{\widetilde{\G}} \def\cS{{\cal S}}

\hfill November 21, 2005 
\begin{center}{\bf On the Limiting Distribution for the Longest Alternating Sequence\\ in a Random Permutation}\end{center}

\begin{center}{\bf Harold Widom}\end{center}

\begin{center}{\it Department of Mathematics\\
University of California\\Santa Cruz, CA 95064\\
e-mail: widom@math.ucsc.edu}\end{center}

\begin{abstract}
Recently Richard Stanley \cite{S} initiated a study of the distribution of the length as$_n(w)$ of the longest alternating subsequence in a random permutation $w$ from the symmetric group $\cS_n$. Among other things he found an explicit formula for the generating function (on $n$ and $k$) for Pr\,(as$_n(w)\le k)$ and conjectured that the distribution, suitably centered and normalized, tended to a Gaussian with variance 8/45. In this note we present a proof of the conjecture based on the generating function.
\end{abstract}

If $w=(w_1\cdots w_n)$ is a permutation in the symmetric group $\cS_n$ then an increasing subsequence of length $k$ is a subsequence $w_{i_1}\cdots w_{i_k}$ satisfying
\[w_{i_1}<w_{i_2}\cdots<w_{i_k}.\]
The random variable is$_n(w)$, the length of the longest increasing sequence in a random permuation $w$ from $\cS_n$, has been much studied. Its mean was first determined asymptotically by Logan-Shepp \cite{LS} and Vershik-Kerov \cite{VK}, proving a conjecture of Ulam, and the limiting distribution was determined in the celebrated work of Beik-Deift-Johansson \cite{BDJ}.

Recently Stanley \cite{S} initiated a study of the distribution of the length of the longest alternating subsequences of $w$, a subsequence such that
\[w_{i_1}>w_{i_2}<w_{i_3}>\cdots w_{i_k}.\]
If as$_n(w)$ denotes the length of the longest alternating sequence in a random permutation $w$ from $\cS_n$ denote by $p_n(k)$ the distribution function for as$_n(w)$,
\[p_n(k)={1\ov n!}\{\#w:{\rm as}_n(w)\le k\}.\]
Stanley found an explicit formula for the generating function 
\[B(x,y)=\sum_{k,n\ge0}p_n(k)\,y^k\,x^n.\]
It is given by
\[B(x,y)={1+\r+2y\,e^{\r x}+(1-\r)\,e^{2\r x}\ov 1+\r-y^2+(1-\r-y^2)\,e^{2\r x}},\]
where $\r=\sqrt{1-y^2}$. Using certain consequences of this formula he made the conjecture that the limit
\[K(t)=\lim_{n\to\iy}p_n(2n/3+\sqrt n\, t).\]
exists. We shall show here that this is so, and that $K(t)$ is a Gaussian with variance $8/45$:
\[K(t)={1\ov\sqrt\pi}\int_{-\iy}^{t\,\sqrt{45}/4}e^{-s^2}\,ds.\]

Two other proofs of the conjecture, very different from this one, are indicated in \cite{S}.

We use the integral representation
\[p_n(k)={1\ov(2\pi i)^2}\int\int B(x,y)\,x^{-n-1}\,y^{-k-1}\,dx\,dy,\]
where the contours of integration are little curves around zero. We begin by making the substitution
\[y=\sech u,\]
where $u$ runs along a vertical line from $a+i\pi$ to $a-i\pi$, with $a$ a large positive real number. Then $y$ runs counterclockwise around a little curve surrounding zero. 

With this substitution 
\[\r=\sqrt {1-\sech^2 u}=\tanh u\]
(when $y$ is small $\r$ is close to 1, which is why $\r$ is given by this square root rather than its negative) and we find that
\[p_n(k)={1\ov(2\pi i)^2}\int\int{e^{2\,u}+2\,e^u\,e^{x\,\tanh u}+e^{2\,x\,\tanh u}\ov e^{2\,u}-e^{2\,x\,\tanh u}}\,\cosh^{k}u\;x^{-n-1}\,dx\,du,\]
where now $u$ runs from $a-\pi i$ to $a+\pi i$.
 
We perform the $x$ integration first. The poles of the integrand are at $x=0$ and
\[x=(u+j\pi i)\,\coth u,\ \ \ \ -\iy<j<\iy.\]
The residue of the quotient in the integrand at the $j$th point is computed to be zero when $j$ is odd (so there is actually no pole there) and $-2\coth u$ when $j$ is even. We evaluate the $x$-integral by integrating around a sequence of expanding contours passing half-way between the poles\footnote{The quotient in the integrand equals 1 plus a constant times $(1+e^{u-x\,\tanh u})/\sinh (u-x\,\tanh u)$. There is an expanding sequence of contours on which this is uniformly bounded. A linear change of variable shows that it is enough to show this for the function $(1+e^z)/\sinh z$. One can take the $j$th contour to be the square with vertices $\pm(j+1/2)\,\pi\pm(j+1/2)\,\pi i$.}
and find that
\be p_n(k)={1\ov \pi i}\int\cosh^{k}u\;\tanh^{n} u\,\sum_{j=-\iy}^\iy\, (u+2j\pi i)^{-n-1}\,du.\label{p}\ee
(Since $j$ was even we replaced it by $2j$.) 
 
Because of the periodicity of the integrand, we may take any real $a$ and as contour of integration any curve from $a-\pi i$ to $a+\pi i$ as long as it stays in the region $|\I u|\le\pi$ and to the right of the poles at zero and $\pm i\pi/2$. 

Here is how we will proceed. Because $u$ lies in the region $|\I u|\le\pi$ it seems likely that the term of the sum corresponding to $j=0$ will dominate the rest. Suppose this is so, and that $k$ will be of the order $n$, say equal to $\a n+o(n)$ for some positive constant $\a$. Then the main part of the resulting integrand will be
\[(\cosh u)^{\a n}\;\left({\tanh u\ov u}\right)^n,\]
which we write as $e^{n\,\s(u)}$. The method of steepest descent, or saddle point method, tells us that we should try to take as our contour of integration one on which $\I \s(u)$ is constant and on which $\R \s(u)$ achieves its maximum at a point $u_0$. If $u_0$ is not an end-point of the curve then $\s'(u_0)=0$. This $u_0$ is a saddle point for the function $|e^{\s(u)}|$ and at all points of the curve its direction away from the saddle point is that of most rapid decrease of the function --- whence the names for the method. The main contribution to the integral will come from the immediate neighborhood of the single saddle point, and so one only has to use an expansion of the integrand near this point to determine the asymptotics. In practice one first locates a saddle point, follows a steepest descent curve until it ends (at a zero of the integrand, at another saddle point, or the boundary of the region which could include $\iy$), and then sees if the original contour may be deformed to it.

In our case the saddle points are the zeros of
\[\s'(u)=\a\,{\sinh u\ov \cosh u}+{1\ov \sinh u\,\cosh u}-{1\ov u}.\]
We shall show that there are three saddle points in the region $|\I u|<\pi$, and that they coincide when $\a=2/3$.

Let us consider $\Delta\,\textrm{arg}\, \s'(u)$ over the boundary of the region, with little semi-circular indentations of diameter $\dl$ below the pole at $i\pi$ and above the pole at $-i\pi$. Using the fact that $\s'(u)$ is odd we need consider only the upper part. We also use the fact that when ${\rm Im}\;u=\pi$ the first two summands in the expression for $\s'(u)$ are real while $-1/u$ has positive imaginary part. As $u$ goes from $+\iy+i\pi$ to $-\iy+i\pi$ along the upper boundary with indentation $\s'(u)$ goes from $\a$ to $\dl\inv+O(1)$ staying in the upper half-plane, then counterclockwise around a large semi-circle (more precisely, within $O(1)$ of one) to $-\dl\inv+O(1)$, then to $-\a$ while again staying in the upper half-plane. Clearly, then,  $\Delta\,\textrm{arg}\, \s'(u)$ over this upper boundary is $\pi$, and so $\Delta\,\textrm{arg}\, \s'(u)$ over the full boundary is $2\pi$. This means that in our region the number of zeros of $\s'(u)$, counting multiplicity, is one more than the number of poles. There are two poles, at $\pm i\pi/2$, so there are three zeros. 

We compute that in the neighborhood of $u=0$ 
\[\s(u)=\a\,\left({1\ov2} u^2-{1\ov12}u^4+O(u^6)\right)-{1\ov3}u^2+{7\ov 90}u^4+O(u^6).\]
When $\a=2/3$ this has a zero of order four at $u=0$, so $\s'(u)$ has a zero of order three and therefore this is the only saddle point. So we do take $\a=2/3$, and near $u=0$
\[\s(u)={1\ov45}u^4+O(u^6).\]

The next step is to determine the steepest descent curve, or at least to describe it in sufficient detail. It is a curve (or curves) emanating from the saddle point on which $\I \s(u)$ is constant (in this case zero) and on which $\R  \s(u)$ decreases as we move away from the saddle point. We see that there are four such curves, emanating from 0 in the directions ${\rm arg}\; u=\pm\pi/4$ in the right half-plane and the directions ${\rm arg}\; u=\pm3\pi/4$ in the left half-plane. By the symmetries of $\s(u)$ (namely that it is even and takes conjugate values at conjugate values of $u$) these curves are reflections of each other in the real and imaginary axes. On all of them $\s(u)$ is real and decreases as we move away from the saddle point. None of the four curves, when extended, can intersect in the interior of the region in question because if they did it would have to be at a saddle point other than $u=0$, and there are none. It follows from this, and the symmetry of the four curves, that the ones emanating from zero in the directions ${\rm arg}\; u=\pm\pi/4$ stay in the right half-plane and meet the upper boundary of our region at points $a\pm i\pi$ with $a\ge0$. The two together, from $a-i\pi$ to 0 to $a+i\pi$, constitute the curve we use as our contour. We denote it by $\G$. Its important properties are its location (to the right of the poles at $\pm i\pi/2$ and therefore giving the same integral as the original contour) and the fact that on it $\s(u)$ is real and decreases from zero as we move away from the saddle point. 

Let us now take $k=2n/3+\sqrt n\, t$ in (\ref{p}). The integrand has a pole at $u=0$, coming from the term of the sum with $j=0$. Since 0 was to the left of the original contour of integration, we replace the term $1/u$ by $1/(u-\ep)$, then deform the contour to $\G$, and then take the limit as $\ep\to0$. This is indicated by replacing the factor $1/u$ by $1/(u-0)$. (The integral arising from this term can also be expressed in terms of a principal value integral.)
Thus (\ref{p}) becomes in present notation
\[p_n(2n/3+\sqrt n\, t)={1\ov \pi i}\int_\G\,e^{n\,\s(u)}\,(\cosh u)^{\sqrt n\, t} \,{du\ov u-0}\]
\[+
{1\ov \pi i}\int_\G\,e^{n\,\s(u)}\,(\cosh u)^{\sqrt n\, t} \,\sum_{j\ne0}\left({u\ov u+2j\pi i}\right)^{n}\,{du\ov u+2j\pi i}.\]

We shall use the notations $c$ and $C$ to denote constants, small and large respectively, which will vary with each use.

We first show that the integral involving the sum is exponentially small. The curve $\G$ lies in some boundeed subregion of $|\I u|\le\pi$. If $|j|$ is sufficiently large then we shall have
\[\left|{u\ov u+2j\pi i}\right|<{1\ov2|j|}\]
uniformly on $\G$, and the sum over these $j$ will be $O(2^{-n})$.
Each of the finitely many other summands has absolute value at most
\[\left|{u\ov u\pm2\pi i}\right|^n,\]
the sign being that of $j$, because  $u\pm2\pi i$ and $u+ 2j\pi i$ have the same real part but the absolute value of the latter is larger than that of the former when $|j|>1$.  The other factor in the integrand, $(u+2j\pi i)\inv$, is bounded. It follows that the absolute value of the integral is at most a constant times
\[\int_\G\,e^{n\,\s(u)}\,|\cosh u|^{\sqrt n\, t}\,
\left(\left|{u\ov u+2\pi i}\right|^n+\left|{u\ov u-2\pi i}\right|^n+2^{-n}\right)\,|du|.\] 

In a sufficiently small neighborhood of zero in $\G$ the expression is parentheses will be $O(e^{-c\,n})$. The cosh factor is $O(e^{C\,\sqrt n})$ and the first factor is at most 1. Thus the product is $O(e^{-c\,n})$. Outside this neighborhood $e^{n\,\s(u)}=O(e^{-c\,n})$ because $\s(u)$ decreases from zero away from the saddle point. The expression in parentheses is $O(1)$ and the cosh factor is as before. Hence the integrand is $O(e^{-c\,n})$ on $\G$ and therefore so is the integral.

We have shown that the integral involving the sum is exponentially small. Let us now consider the other integral.

Take any $\dl>0$ and denote by $\G_1$ the part of $\G$ on which $|u|<n^{-\dl}$ and by $\G_2$ the remainder of $\G$. Then
$e^{n\,\s(u)}=O(e^{-c\,n^{1-4\dl}})$ on $\G_2$, for this is the estimate when $|u|=n^{-\dl}$ and $e^{n\,\s(u)}$ is less than this on $\G_2$. (Notice that ${\rm arg}\,u^4$ is close to $\pm\pi$  when $|u|=n^{-\dl}$.) Including the other factors in the integrand, we see that on $\G_2$ the integrand is $O(e^{-c\,n^{1-4\dl}+C\,n^{1/2}}\,n^\dl)$. If we choose $\dl<1/8$, and we do, then this bound is just $O(e^{-c\,n^{1-4\dl}})$. Hence the integral over $\G_2$ is exponentially small. 

Finally we come to the part that gives the asymptotics, the integral over $\G_1$. We use the expansion
\[e^{n\,\s(u)}\,(\cosh u)^{\sqrt n\, t}=\exp\Big\{n\,(u^4/45+O(u^6))+\sqrt n\,t\,(u^2/2+O(u^4))\Big\}\]
valid near $u=0$.
We want to replace $\G_1$ by $\tlG$ which consists of the portion of the rays ${\rm arg}\,u=\pm\pi/4$ from zero to the points where $|u|=n^{-\dl}$. The difference between the two integrals is the integral over the line segment joining the ends of $\G_1$ and $\tlG$, on which the integrand is exponentially small. So, with exponentially small error, we may make this replacement. The contribution of the integral over $\tlG$ may be written
\[{1\ov \pi i}\int_{\tlG}\,\exp\Big\{n\,(u^4/45+O(u^6))+\sqrt n\, t\,(u^2/2+O(u^4))\Big\}{du\ov u-0}.\]
The change of variable $u=n^{-1/4}\,v^{1/2}$ replaces this by
\[{1\ov 2\pi i}\int_{-i\,n^{1/2-2\dl}}^{i\,n^{1/2-2\dl}}\,\exp\Big\{(v^2/45+O(n^{-1/2}\,v^3))+ t\,(v/2+O(n^{-1/2}\,v^2))\Big\}{dv\ov v-0}\]
\[={1\ov 2\pi i}\int_{-i\,n^{1/2-2\dl}}^{i\,n^{1/2-2\dl}}\,\exp\Big\{(v^2/45+O(n^{-2\dl}\,v^2))+ t\,(v/2+O(n^{-2\dl}\,v))\Big\}{dv\ov v-0}.\]

The exponential factor is uniformly bounded by an integrable function $e^{-c\,|v|^2+C\,|v|}$ so\footnote{This is not quite enough. An integral $\int_{-i\iy}^{i\iy} f_n(v)\,dv/(v-0)$ is equal to 
\[\pi i\,f_n(0)+{\rm PV}\int_{-i\iy}^{i\iy} f_n(v)\,dv/v
=\pi i\,f_n(0)+\int_{0}^{i\iy} (f_n(v)-f_n(-v))\,dv/v.\]
To show convergence to the limit $\int_{-i\iy}^{i\iy} f(v)\,dv/(v-0)$ it is enough to have pointwise convergene of $f_n$ to $f$, dominated convergence outside a neighborhood of 0 and, for example, uniform boundedness of the derivatives $f_n'$ in this neighborhood. This is seen to hold in our case.} we can take the limit under the integral and find that the above has limit
\[{1\ov 2\pi i}\int_{-i\,\iy}^{i\,\iy}\,e^{v^2/45+ t\,v/2}{dv\ov v-0}.\]
This is equal to
\[{1\ov\sqrt\pi}\int_{-\iy}^{t\,\sqrt{45}/4}e^{-s^2}\,ds.\]

\begin{center}{\bf Acknowlegments}\end{center}
The author thanks Richard Stanley for making his paper \cite{S} available to the author before publication, and Craig Tracy for alerting the author to the question. Research was supported by the National Science Foundation under grant DMS-0243982.

\end{document}